\documentclass[11pt]{amsart}

\usepackage[utf8]{inputenc}
\usepackage[T1]{fontenc,url}

\usepackage[left=3.3cm,right=3.3cm,bottom=2.3cm,top=2.3cm]{geometry}

\usepackage{amsmath,amsthm,amsfonts,amssymb}
\usepackage{mathrsfs}
\usepackage[shortlabels]{enumitem}
\usepackage{cite}
\usepackage{mathtools}
\usepackage{color}
\definecolor{shadecolor}{gray}{0.875}
\definecolor{col}{RGB}{42, 95, 151}
\usepackage[colorlinks=true, linkcolor=col, citecolor=col, filecolor=col, menucolor=col, urlcolor=col]{hyperref}

\numberwithin{equation}{section}

\theoremstyle{plain}
\newtheorem{theorem}{Theorem}[section]
\newtheorem*{lemma*}{Lemma}

\newtheorem*{theorem*}{Theorem}
\newtheorem{proposition}[theorem]{Proposition}
\newtheorem*{proposition*}{Proposition}

\newtheorem*{corollary*}{Corollary}

\theoremstyle{definition}
\newtheorem{remark}[theorem]{Remark}
\newtheorem*{remark*}{Remark}
\newtheorem*{definition*}{Definition}

\newtheorem*{example*}{Example}

\newtheorem*{question*}{Question}

\def\c1{\operatorname{c_1}}
\def\c2{\operatorname{c_2}}

\def\CC{{\mathbb C}}
\def\ZZ{{\mathbb Z}}

\def\PP{{\mathbb P}}

\def\O{{\mathcal O}}

\def\+{\oplus}
\def\*{\otimes}

\newcommand{\CH}{\operatorname{CH}}

\def\PP{{\mathbb P}}
\def\ZZ{{\mathbb Z}}

\def\CC{{\mathbb C}}
\def\FF{{\mathbb F}}
\def\O{{\mathcal O}}

\title[Cohomological and Chow decompositions of the diagonal]{A remark on cohomological and Chow decompositions of the diagonal}
\author{John Christian Ottem}
 \address{Department of Mathematics,
 University of Oslo,
 Box 1053, Blindern,
 0316 Oslo, Norway}
 \email{johnco@math.uio.no}

\begin{document}

\begin{abstract}
We show that there are quartic double solids whose desingularizations admit a cohomological decomposition of the diagonal but no Chow decomposition.
\end{abstract}
\maketitle

\section{Introduction}
Let $X$ be a smooth projective variety of dimension $n$ over an algebraically closed field $k$. 
We say that $X$ admits a \emph{Chow decomposition of the diagonal} if the class of the diagonal $\Delta_X\subset X\times X$ in the Chow group $\CH^n(X\times X)$ can be written as \begin{equation}\label{eq:Chow}
 [\Delta_X]=[X\times p]+[\Gamma],
\end{equation}
where $p\in X(k)$ and $\Gamma$ is a cycle supported on $D\times X$ for a closed subset $D\subsetneq X$.
This is equivalent to \emph{universal $\CH_0$-triviality}, meaning that the degree map
$\CH_0(X_L)\to\ZZ$ is an isomorphism for every field extension $L/k$.

Similarly, if $k=\CC$, a \emph{cohomological decomposition of the diagonal} is the same equality \eqref{eq:Chow}, but now in $H^{2n}(X\times X,\ZZ)$. 
If $X$ has a Chow decomposition, then it has a cohomological decomposition, by the cycle class map. 

If $X$ is rational, or more generally stably rational, then it admits a Chow decomposition of the diagonal. 
In \cite{Voi}, Voisin introduced the specialization method using the failure of such a decomposition as an obstruction to stable rationality.
This method was later extended by Colliot-Th\'el\`ene--Pirutka
\cite{CTP} and Schreieder \cite{SchDuke,Sch}, leading to many results on irrationality of hypersurfaces in projective space \cite{CTP,EGFS,Tot,Sch,LS,Ott,Hot,BBPZ}.
In these applications, one usually rules out a Chow decomposition, which is a priori a weaker conclusion than ruling out a cohomological
decomposition. One reason is that a Chow decomposition can be obstructed by studying the action of the diagonal on other cohomology theories, such as unramified cohomology. 

Voisin asked whether a smooth complex projective variety can admit a cohomological decomposition without a Chow 
decomposition \cite[Question~3.27]{VoiSurvey}. In this note, we answer this question affirmatively.

\begin{theorem}\label{thm:main}
There exists a quartic double solid with 10 nodes whose desingularization admits a cohomological decomposition of the diagonal but no Chow decomposition.
\end{theorem}
So the two notions are different even for rationally connected threefolds. For comparison, they
are equivalent for smooth cubic threefolds and fourfolds \cite[Theorem~1.1]{VoiCubic}. Cohomological decompositions were used recently by Voisin \cite[Theorem~1.3]{VoiFano} to show that, for several classes of Fano threefolds, the stable birational type of a very general member determines it up to isomorphism.

Our example is in fact closely related to Voisin's original application of the specialization method \cite{Voi}. In that paper, Voisin shows that a very general quartic double solid is not stably rational by specializing to a double solid of Artin--Mumford type. This
threefold has 10 nodes and its desingularization has nontrivial $2$-torsion in $H^3$, which obstructs having a decomposition of the diagonal. So although a smooth quartic double solid has torsion free $H^3$, the torsion in the special fiber obstructs a Chow decomposition of the diagonal, by specialization.

The construction in this paper also uses 10-nodal quartic double solids, but we instead arrange that the desingularization $X$ has $H^3(X,\ZZ)=0$. 
Such threefolds automatically have a cohomological decomposition of the diagonal (Proposition~\ref{prop:cohomological}).
To rule out a Chow decomposition, we apply the method of Koll\'ar--Totaro \cite{Kol,Tot}, using reduction to characteristic $2$.
In characteristic $2$, the double cover is purely inseparable and its desingularization admits a nonzero global $2$-form, which obstructs a Chow decomposition. Using the specialization argument of \cite{Tot}, we see that the complex threefold does not admit a Chow decomposition either.

\smallskip

\noindent{\bf Acknowledgements}. I thank Claire Voisin for many helpful discussions related to the topics of this paper, and in particular for explaining the consequences mentioned in Section~\ref{sec:other-consequences}.

\section{Quartic double solids with 10 nodes}\label{ianrinniasdfnvar}
The following proposition gives the main strategy for finding threefolds with a cohomological decomposition of the diagonal.

\begin{proposition}\label{prop:cohomological}
Let $X$ be a smooth projective rationally connected complex threefold. If $H^3(X,\ZZ)=0$, then $X$ admits a cohomological decomposition of the diagonal.
\end{proposition}
\begin{proof}
Since $X$ is rationally connected, it is simply connected, so $H^1(X,\ZZ)=0$ and $H^2(X,\ZZ)$ is torsion free.
Therefore $H^\ell(X,\ZZ)$  is torsion free for all $\ell$ and $H^\ell(X,\ZZ)=0$ for all odd $\ell$.
Moreover, the group $H^2(X,\ZZ)$ is generated by classes of divisors, and $H^4(X,\ZZ)$ is generated by classes of curves, by Voisin's result on the integral Hodge conjecture \cite[Theorem~2]{VoiHodge}. 
Let $D_i$ be divisors forming a basis of $H^2(X,\ZZ)$ and $C_i$ be 1-cycles forming a dual basis for $H^4(X,\ZZ)$.
The K\"unneth formula now shows that
\begin{equation}\label{aidvnrisarnnnfjwefnui}
 [\Delta_X]=[X\times p]+[p\times X]
  +\sum_i\bigl([D_i\times C_i]+[C_i\times D_i]\bigr).
\end{equation}
So $X$ admits a cohomological decomposition of the diagonal.
\end{proof}
To apply this, we consider desingularizations of nodal quartic double solids.

Let $F\in\CC[x_0,\ldots,x_3]$ be a homogeneous polynomial of degree $4$ and let $S\subset\PP^3_\CC$ be the corresponding quartic surface.
We assume that the only singularities of $S$ are 10 nodes $\Sigma\subset S$. 
Then the double cover $W\to\PP^3_\CC$, given by $w^2=F$, also has 10 nodes.
Let $X\to W$ be the blow-up of these nodes. Then $X$ is a smooth projective rationally connected threefold.

The third Betti number $b_3(X)$ is computed using the {\em defect}, which in this case measures the 
failure of the nodes to impose independent conditions on quadrics \cite[3.16]{Cl}. 
More precisely, for 10 nodes, we have $b_3(X)=2h^0(\PP^3,I_\Sigma(2))$, so $b_3(X)=0$ precisely when the 10 nodes are not contained in a quadric.
Under this assumption, it remains to check that $H^3(X,\ZZ)$ has no nonzero torsion.

The torsion in $H^3(X,\ZZ)$ is closely related to the geometry of the quartic surface $S$ and its nodes.  
If we project from a node, we obtain a generically $2$-to-$1$ rational map $S\dashrightarrow\PP^2_\CC$, branched along a plane curve  
$B\subset \PP^2_\CC$ of degree 6. The curve $B$ is also the discriminant of the conic bundle $W\dashrightarrow\PP^2_\CC$.

Suppose that $B$ is irreducible with exactly 9 nodes and that the conics above these nodes have rank 2. 
Blowing up $\PP^2_\CC$ along these nodes, and making a birational modification of the total space, we obtain
a standard conic bundle $X'\to \operatorname{Bl}_9 \PP^2_\CC$ (see \cite[§4, Theorem~2 and p.~94]{AM}). 
The new discriminant is the strict transform of $B$, which is smooth and irreducible.
The formula for the Brauer group of a standard conic bundle \cite[Theorem~11.3.14]{CTS} therefore gives
$\operatorname{Br}(X')=0$ and hence $\operatorname{Br}(X)=0$. Since $X$ is rationally connected, $\operatorname{Br}(X)\simeq H^3(X,\ZZ)_{\mathrm{tors}}$.
Therefore, if $B$ satisfies these conditions and the 10 nodes do not lie on a quadric, then $H^3(X,\ZZ)=0$.

\begin{remark}
The study of nodal quartic surfaces goes back to Cayley \cite[p.~198]{Cay}. Cayley distinguished between
10-nodal quartics whose projection discriminant is an irreducible curve of degree 6 and those for which it splits into two cubic curves. The second case corresponds to 
\emph{quartic symmetroids}, which are defined by the determinant of a symmetric $4\times4$ matrix of linear forms. 
The corresponding quartic double solids were used by Artin--Mumford \cite{AM} in their examples of unirational threefolds with nonzero torsion in $H^3$.
\end{remark}
To rule out a Chow decomposition, we apply the method of Koll\'ar--Totaro \cite{Kol,Tot} using reduction to characteristic $2$. 

Suppose that $F$ has coefficients in $\ZZ$, and write $f$ for its reduction to $k=\overline{\FF}_2$. 
We assume that the critical points of $f$ are \emph{almost nondegenerate} in the sense of Koll\'ar \cite{Kol}.
Explicitly, write $f$ in local coordinates near a critical point $p$ as $f=f(p)+q_2+q_3+\cdots$ with each $q_i$ 
homogeneous of degree $i$. Then the condition is that the Hessian of $q_2$ has rank $2$ and $q_3$ is nonzero on its 1-dimensional kernel. 

Let $\pi\colon Y\to\PP^3_k$ be the inseparable double cover $w^2=f$. 
Blowing up the singular points of $Y$ gives a resolution $\varphi\colon Z\to Y$. Its exceptional divisors are quadric cones, which are universally $\CH_0$-trivial since $k$ is algebraically closed. Hence the morphism $\varphi$ is universally $\CH_0$-trivial as well. Koll\'ar's construction \cite[\S\S3.1 and~4.2]{OkaCyc} now gives an injection $\O_Z\hookrightarrow\Omega_Z^2$. Hence $H^0(Z,\Omega_Z^2)\ne0$, so $Z$ is not universally $\CH_0$-trivial
\cite[Lemma~2.2]{Tot}. Now, the double cover $w^2=F$ of $\PP^3_{\ZZ_{(2)}}$ is flat and projective with geometrically integral fibers.
Since $Z$ is not universally $\CH_0$-trivial, applying the specialization argument \cite[Theorem~2.3]{Tot}, after base change to the ring of Witt vectors $W(k)$,  
we conclude that $X$ admits no Chow decomposition of the diagonal.
 
%

\subsection{Explicit equations}\label{subsec:explicit-equations}
For the main example, we look for a quartic polynomial $F$ with the following properties:
\begin{enumerate}
\renewcommand{\labelenumi}{(\theenumi)}
\item $S$ has exactly 10 ordinary nodes, imposing independent conditions on quadrics;
\item projection from one of the nodes of $S$ gives an irreducible 9-nodal discriminant, with rank-2 conics over its nodes;
\item the reduction $f$ of $F$ modulo 2 has only almost nondegenerate critical points.
\end{enumerate}
For convenience, we view $\PP^3$ as a linear subspace of $\PP^4$. 
Let $e_j$ denote the $j$th elementary symmetric polynomial in $x_0,\ldots,x_4$. We start in characteristic $2$ with the quartic surface
$$
 S_0=\{e_1=e_4=0\}\subset\PP^4_k.
$$
Write $H_0=\{e_1=0\}\simeq\PP^3_k$ and $f=e_4|_{H_0}$.
A Jacobian computation shows that the critical points of $f$ are the 10 permutations of $(0:0:0:1:1)$.

To check the almost nondegeneracy  condition, it suffices to consider the point $(0:0:0:1:1)$. On the affine chart $x_3=1$, we eliminate $x_4$ using $x_4=1+x_0+x_1+x_2$. 
Then $f=q_2+q_3+q_4$ where
$$
 q_2=x_0x_1+x_0x_2+x_1x_2, \qquad q_3=(x_0+x_1+x_2)q_2.
$$
The Hessian of $q_2$ has rank 2 with kernel generated by $(1,1,1)$, and $q_3(t,t,t)=t^3$, so we see that the critical points of $f$ are almost nondegenerate.

We now define $F$ by taking the following lift of $H_0$ and $f$:
\begin{equation}\label{eq:construction}
 H=\{x_0+x_1+x_2+3x_3+3x_4=0\}\subset\PP^4_\ZZ, \qquad F=e_4+4x_3^2x_4^2.
\end{equation}
One checks that the singular locus of the surface $S=\{F=0\}\subset H_\CC$ consists of 10 nodes, obtained from 
$(1:-1:0:0:0)$, $(3:0:0:-1:0)$ and $(2:2:2:-1:-1)$ by permuting $x_0,x_1,x_2$ and switching $x_3,x_4$. These 10 points reduce modulo 2 to the 10 points above and they are not contained in any quadric surface in $\PP^3_\CC$.

Finally, projecting from the node $(2:2:2:-1:-1)$, in coordinates $(u:v:w)=(x_0+x_3+x_4:x_1+x_3+x_4:x_3-x_4)$, the discriminant curve $B$ is defined by
\begin{equation}\label{eq:sextic}
 w^2(3u^2+3uv+3v^2-5w^2)(u^2+uv+v^2+w^2)-4u^2v^2(u+v)^2=0.
\end{equation}
One computes that this curve is irreducible with exactly 9 ordinary nodes, and the conics over these nodes all have rank 2. Therefore, (1)--(3) hold and Theorem~\ref{thm:main} follows.

\section{Other consequences}\label{sec:other-consequences}

Claire Voisin pointed out the following two consequences of Theorem~\ref{thm:main}. In the statement below, $\circ$ denotes composition of correspondences.

\begin{proposition}\label{prop:other-consequences}
Let $X$ be the threefold of Theorem~\ref{thm:main}. Then:
\begin{enumerate}
\item there exists a nonzero correspondence $E\in\CH^3(X\times X)$
such that $[E]=0$ in $H^6(X\times X,\ZZ)$ and $E^{\circ2}=E$;
\item $X$ admits no decomposition of the diagonal modulo algebraic equivalence.
\end{enumerate}
\end{proposition}

\begin{proof}
We follow Voisin's argument in
\cite[proof of Proposition~2.1]{VoiCubic} using the decomposition \eqref{aidvnrisarnnnfjwefnui}. 
Let $D_i$ and $C_i$ be cycles as in Proposition \ref{prop:cohomological}, so that $\deg(D_i\cdot C_j)=\delta_{ij}$.  For a point $p\in X$, let $E=\Delta_X-X\times p-\Gamma$ where
$$
\Gamma=p\times X+\sum_i(D_i\times C_i+C_i\times D_i).
$$
By the intersection numbers, $\Gamma^{\circ2}=\Gamma$, and $\Gamma\circ(X\times p)=(X\times p)\circ\Gamma=0$.
Since $X\times p$ is also idempotent, $E^{\circ2}=E$. 
By \eqref{aidvnrisarnnnfjwefnui}, $[E]=0$. Moreover, $E\neq 0$, because $X$ has no Chow decomposition of the diagonal. This proves (1).

Statement (2) follows from \cite[Proposition~2.1]{VoiCubic}, which states that any decomposition of the diagonal modulo algebraic equivalence lifts to a Chow decomposition.
This result uses the nilpotence of algebraically trivial correspondences, due to Voisin \cite{VoiNil} and Voevodsky \cite{Voe}.
\end{proof}

Part (1) shows that nilpotence for homologically trivial correspondences fails with integral coefficients, since in fact $E^{\circ N}=E\ne0$ for every $N\ge1$.
Soul\'e and Voisin had already constructed homologically trivial integral cycles which are not nilpotent for the smash product \cite[Theorem~5]{SV}. 

\bigskip

A \verb|Macaulay2| file verifying some of the calculations in this paper is available at
\begin{center}
\url{https://www.mn.uio.no/math/personer/vit/johnco/papers/10nodes.m2}
\end{center}

\smallskip

\noindent {\bf AI disclosure}. ChatGPT-6 was used to assist with finding the explicit example, checking computations and revising the paper.


\begin{thebibliography}{99}

\bibitem{AM}
M. Artin and D. Mumford,
\emph{Some elementary examples of unirational varieties which are not rational},
Proc. London Math. Soc. (3) \textbf{25} (1972), 75--95.

\bibitem{BBPZ}
C. B\"ohning, H.-C. Graf von Bothmer, A. Pirutka and Z. Zhang,
\emph{Stable irrationality of quartic hypersurfaces},
preprint (2026),
\href{https://arxiv.org/abs/2611.15702}{arXiv:2611.15702}.

\bibitem{Cay}
A. Cayley,
\emph{A second memoir on quartic surfaces},
Proc. London Math. Soc. (1) \textbf{3} (1869--1871), 198--202.

\bibitem{Cl}
C. H. Clemens,
\emph{Double solids},
Adv. in Math. \textbf{47} (1983), 107--230.

\bibitem{CTP}
J.-L. Colliot-Th\'el\`ene and A. Pirutka,
\emph{Hypersurfaces quartiques de dimension 3 : non-rationalit\'e stable},
Ann. Sci. \'Ec. Norm. Sup\'er. (4) \textbf{49} (2016), 371--397.

\bibitem{CTS}
J.-L. Colliot-Th\'el\`ene and A. N. Skorobogatov,
\emph{The Brauer--Grothendieck Group},
Ergebnisse der Mathematik und ihrer Grenzgebiete (3), vol.~71,
Springer, Cham, 2021.


\bibitem{EGFS}
P. Engel, O. de Gaay Fortman and S. Schreieder,
\emph{Matroids and the integral Hodge conjecture for abelian varieties},
preprint (2025),
\href{https://arxiv.org/abs/2507.15704}{arXiv:2507.15704}.

\bibitem{Hot}
J. Hotchkiss,
\emph{The rationality problem for hypersurfaces of degree at least five},
preprint (2026),
\href{https://arxiv.org/abs/2610.10509}{arXiv:2610.10509}.

\bibitem{Kol}
J. Koll\'ar,
\emph{Nonrational hypersurfaces},
J. Amer. Math. Soc. \textbf{8} (1995), no.~1, 241--249.

\bibitem{LS}
J. Lange and S. Schreieder,
\emph{On the rationality problem for low degree hypersurfaces},
Forum Math. Pi \textbf{13} (2025), e25.

\bibitem{OkaCyc}
T. Okada,
\emph{Stable rationality of cyclic covers of projective spaces},
Proc. Edinb. Math. Soc. (2) \textbf{62} (2019), 667--682.

\bibitem{Ott}
J. C. Ottem,
\emph{Very general quartic sixfolds are not stably rational},
preprint (2026),
\href{https://arxiv.org/abs/2609.10231}{arXiv:2609.10231}.

\bibitem{SchDuke}
S. Schreieder,
\emph{On the rationality problem for quadric bundles},
Duke Math. J. \textbf{168} (2019), no.~2, 187--223.

\bibitem{Sch}
S. Schreieder,
\emph{Stably irrational hypersurfaces of small slopes},
J. Amer. Math. Soc. \textbf{32} (2019), 1171--1199.

\bibitem{SV}
C. Soul\'e and C. Voisin,
\emph{Torsion cohomology classes and algebraic cycles on complex projective manifolds},
Adv. Math. \textbf{198} (2005), 107--127.

\bibitem{Tot}
B. Totaro,
\emph{Hypersurfaces that are not stably rational},
J. Amer. Math. Soc. \textbf{29} (2016), 883--891.

\bibitem{VoiHodge}
C. Voisin,
\emph{On integral Hodge classes on uniruled or Calabi--Yau threefolds},
in \emph{Moduli Spaces and Arithmetic Geometry},
Advanced Studies in Pure Mathematics, vol.~45,
Mathematical Society of Japan, Tokyo, 2006, 43--73.

\bibitem{Voe}
V. Voevodsky,
\emph{A nilpotence theorem for cycles algebraically equivalent to zero},
Internat. Math. Res. Notices \textbf{1995} (1995), no.~4, 187--198.

\bibitem{VoiNil}
C. Voisin,
\emph{Remarks on zero-cycles of self-products of varieties},
in \emph{Moduli of Vector Bundles} (Sanda, 1994; Kyoto, 1994),
Lecture Notes in Pure and Applied Mathematics, vol.~179,
Marcel Dekker, New York, 1996, 265--285.

\bibitem{Voi}
C. Voisin,
\emph{Unirational threefolds with no universal codimension $2$ cycle},
Invent. Math. \textbf{201} (2015), 207--237.

\bibitem{VoiCubic}
C. Voisin,
\emph{On the universal $\CH_0$ group of cubic hypersurfaces},
J. Eur. Math. Soc. \textbf{19} (2017), 1619--1653.

\bibitem{VoiFano}
C. Voisin,
\emph{On the stable birational type of Fano threefolds},
preprint (2026),
\href{https://arxiv.org/abs/2611.14417}{arXiv:2611.14417}.

\bibitem{VoiSurvey}
C. Voisin,
\emph{Birational invariants and decomposition of the diagonal},
in \emph{Birational Geometry of Hypersurfaces},
A. Hochenegger, M. Lehn and P. Stellari (eds.),
Lecture Notes of the Unione Matematica Italiana, vol.~26,
Springer, Cham, 2019, 3--71.

\end{thebibliography}
\end{document}